\documentclass[12pt]{article}
\usepackage{index} %\usepackage{russcorr}
\usepackage[cp1251]{inputenc}
\usepackage[russian]{babel}
\usepackage{amssymb,latexsym,amsmath}
\usepackage{oldlfont,euscript}
\usepackage{amsfonts}
\usepackage{enumerate}
\usepackage{array,longtable}
\usepackage{pstricks-add,xcolor}

\newcommand{\D}{\displaystyle}
\makeindex
\renewindex{default}{idx}{ind}{\hspace{4.0cm}Предметный указатель\label{I}}

\oddsidemargin\evensidemargin    \parindent 0.9cm
\begin{document}\noindent 47A52; 65F22; 65Q30

\vspace{0.3cm}\centerline{\sc  O.V. Matysik, P.P. Zabreiko}

\vspace{0.3cm}\centerline{\bf M.A. KRASNOSELSKII THEOREM}
\centerline{\bf AND ITERATION METHODS OF SOLVING}
\centerline{\bf ILL-POSED \ \ LINEAR \ \ PROBLEMS}
\centerline{\bf WITH A SELF-ADJOINT OPERATOR}

\vspace{0.5cm}

{\small {\bf Abstract:} The article deals with iterative methods of solving linear operator equations $x = Bx + f$ and $Ax = f$ with self-adjoint operators in Hilbert space $X$ in critical case when $\rho(B) = 1$ and $0 \in {\rm Sp}\, A$. The main results are based on the use of M.A. Krasnosel'ski\u{\i} theorem about the convergence of the successive approximations and some its modifications and refinements.}

\vspace{0.5cm}

The method of successive approximations is one of the main methods of approximate solution of linear operator equations of the second order $x = Bx + f$ in Hilbert and Banach spaces. The main theorems on the convergence of this method, the convergence rate, error estimates, etc. is reduced to studying the properties of Neumann series $\D\sum_{n=0}^\infty B^n$ for the corresponding operator $B$ и and are expounded in numerous textbooks and monographs, among which we can mention here \cite{KanAki,MAK...}. At the same time the greater part of the obtained results are related to the so-called noncritical case, when spectral radius $\rho(B)$ of this linear operator is strictly less than $1$, and this condition is necessary and sufficient for the convergence of Neumann series in the space of operators. However, it was later found out that Neumann series can also converge (not at the norm of operators, but only strongly) in the cases when spectral radius $\rho(B)$ of the corresponding operator is equal to $1$. One of the first results in this direction was obtained by M.A. Krasnosel'skii \cite{MAK} (See also \cite{MAK...}), who showed that for equation $x = Bx + f$ with a self-adjoint operator $B$ in Hilbert space in case of $\rho(B) = \|B\| = 1$ and the supplementary hypothesis that $-1$ is not the eigenvalue of $B$ the successive approximations converge to one of the solutions to the equation under study, only if the equation is solvable. This theorem is not trivial, for under the suggested assumptions equation $x = Bx + f$, is, generally speaking, ill-posed.

The {\it primary} aim of the present paper is to show that the theorem of M.A. Krasnosel'skii on the successive approximations convergence for the equations with self-adjoint operators mentioned above contains in itself, with some natural additions, the main results of iteration methods of approximate solving ill-posed linear equations of the second order with self-adjoint operator $B$ in Hilbert space $X$.

The method of successive approximations is also widely used for approximate constructing the solutions to operator equations of the first order $Ax = y$. Here the main scheme of using the successive approximation method is based on the transition from the original equation $Ax = y$ to the equivalent (or almost equivalent) equation of the second order рода $x = Bx + f$. One of the main methods of such transition is based on the use of functions $f(A)$ of operator $A$ as operator $B$. This type of methods has been studied by many authors (See, for instance, \cite{MAK...}), however, the main results here are obtained for correct equations $Ax = y$, to put in other words, under the additional assertion that $0 \notin {\rm Sp}\, A$.

The {\it second} aim of the present paper is to study the possibility of using the method of successive approximations for finding the approximate solutions to equation $Ax = y$, namely, when this equation is incorrect, that is, when $0 \in {\rm Sp}\, A$. We shall confine ourselves to the case when operator $A$, as well as operator $B = f(A)$ constructed on its basis, are self-adjoint operators in Hilbert space $X$. Thus, the main research tool for ill-posed linear operator equations of the first order will be again the above-mentioned theorem of M.A. Krasnosel'skii.

It should be noted that the the theory of ill-posed linear operator equations with a self-adjoint operator in a Hilbert space was developed independently from different aspects and with a different degree of exactness by various authors. It is sufficient to mention the following monographs
\cite{IvaVasTan,TihArs,Lav,VaiVer,BakGon,SamVab,GilGol,SavMat,Mat}. The {\it third} aim of our paper is to show how the results of these works known before make up a general scheme, as well as to formulate a range of new propositions.

We shall list the main problems examined in this paper for successive approaches to the exact solutions to the equations under study. They concern the convergence conditions of these successive approximations and determining their convergence rate for exact solutions in the original and "weakened" \ norms both in the whole space and in some subspaces densely enclosed in the original one; the analysis of the behaviour of residuals and corrections in constructing these successive approximations; and, finally, the behaviour of the corresponding errors in cases when the right-hand members are defined approximately, and when the calculations are made with some mistakes.

\vspace{0.5cm}

\centerline{\bf \S\ 1. The equalities of the second order}

\vspace{0.5cm}

{\bf 1.1. The convergence of successive approximations.} Let $X$ be Hilbert space, and $B$ be a self-adjoint operator, $f \in X$. Consider the equation
 \begin{equation}\label{101}
 x = Bx + f.
 \end{equation}
To find  solutions  to this equation, it is natural to use the method of successive approximations
 \begin{equation}\label{102}
 x_{n+1} = Bx_n + f \qquad (x_0 \in X, \ n = 0,1,2,\ldots).
 \end{equation}
Actually, if the sequence $(x_n)$ defined by the equation (\ref{102}) is convergent, then its limit will be the solution of the equation (\ref{101}).

The convergence analysis of successive approximations (\ref{102}) is carried out with sufficient completeness in case of $\rho(B) < 1$. The latter inequality (for any continuous linear and not obligatorily self-adjoint operator $B$) is equivalent to Neumann series convergence and to equation
 \begin{equation*}
 (I - B)^{-1} = \sum_{n=0}^\infty B^n.
 \end{equation*}
That is why we will further be interested only in the <<critical>> case, when $\rho(B) = 1$. For the self-adjoint operator it means that $\|B\| = 1$ и что ${\rm Sp}\, B \cap \{-1,1\} \ne \emptyset$. In case of   $1 \notin {\rm Sp}\, B$ equation (\ref{101}) remains uniquely solvable at any $f \in X$, though the issue of convergence to the corresponding solution of successive approximations \ref{102}) remains open-ended. At the same time, in case of $1 \in {\rm Sp}\, B$ equation (\ref{101}) proves to be solvable for some right-hand members $f \in X$ (generally speaking, in this case the solution is non-unique) and unsolvable with other right-hand members $f \in X$.

In the cited paper \cite{MAK} M.А.Krasnosel'skii gave an exhaustive answer concerning the conditions of successive approximations convergence in the critical case described above. In the modification of M.A. Krasnoselskii theorem to be given further we can find the statement as to which of solutions the successive approximations converge when equation (\ref{101}) has an ambiguous solution.

{\bf Theorem 1.1.} {\it Let $B$ be a self-adjoint operator with $\rho(B) = 1$ in Hilbert space $X$, while  $-1$ is not its eigenvalue.

Let equality {\rm (\ref{101})} be solved. Then successive approximations {\rm (\ref{102})} at any initial condition $x_0 \in X$ converge to one of the solutions of equation {\rm (\ref{101})}.

More exactly, approximations {\rm (\ref{102})} сonverge to the solution $x_*$ to equation {\rm (\ref{101})}, for which $Px_* = Px_0$, where $P$ is an orthoprojection on the set of eigenvectors of operator $B$, with the eigenvalue $1$.}

$\square$ We will give a simple proof scheme of this theorem (сompare with the one given in \cite{MAK,MAK...}). From (\ref{101}) and  (\ref{102}) there obviously follow equation
 \begin{equation}\label{103}
 x_n = B^nx_0 + (E + B + \ldots + B^{n-1})f \qquad (n = 0,1,2,\ldots),
 \end{equation}
 \begin{equation}\label{104}
 x_* = B^nx_* + (E + B + \ldots + B^{n-1})f \qquad (n = 0,1,2,\ldots),
 \end{equation}
resulting in
 \begin{equation}\label{105}
 x_n - x_* = B^n(x_0 - x_*),
 \end{equation}
Hence,  by virtue of the theorems on spectral decomposition of of self-adjoint operators in Hilbert space (See, for instance, \cite{KanAki,DanShv}),
 \begin{equation}\label{106}
 \|x_n - x_*\|^2 = \int\limits_{{\rm sp}\, B} |\lambda|^{2n} \, (dE_\lambda(x_0 - x_*),x_0 - x_*)),
 \end{equation}
where $E_\lambda$ is a spectral measure for operator $B$. The sequence $|\lambda|^{2n}$ сonverges to zero everywhere on $(-1,1) \cap \, {\rm sp}\, B$. Point $-1$ (if it is enclosed in ${\rm sp}\, B$) under Theorem 1.1, has a zero spectral measure. Point $1$ (if it is again in ${\rm sp}\, B$) can have positive measure, yet only when $Px_0 \ne Px_*$. Тhus,the proposition of Theorem 1.1 results from Lebesgue theorem on the passage to the limit  under integral sign. $\blacksquare$

We observe that the convergence of successive approximations can generally be аrbitrarily poor. It becomes quite obvious from equation (\ref{106}). The corresponding examples can be easily given.

\vspace{0.3cm}

{\bf 1.2. The convergence of residuals and corrections.} We shall now consider the behaviour of residuals $x_n - Bx_n - f$ for approximations (\ref{102}). It follows,
 \begin{equation*}
 x_n - Bx_n - f = x_n - x_{n+1},
 \end{equation*}
namely, the residuals in the case under consideration coincide with the corrections taken with the reversed sign. From (\ref{103}) it follows
 \begin{equation}\label{107}      x_n - Bx_n - f = B^n(x_0 - Bx_0 - f).
 \end{equation}
Once again, by virtue of the spectral theorem for self-adjoint operators this equation brings about equations
 \begin{equation}\label{108}
 \|x_n - Bx_n - f\|^2 = \int\limits_{{\rm sp}\, B} |\lambda|^{2n} \, (dE_\lambda(x_0 - Bx_0 - f),x_0 - Bx_0 - f).
 \end{equation}
We can once again apply Lebesgue theorems of the limiting process to this equation. As a result, we obtain the following assertion:

{\bf Theorem 1.2.} {\it Let $B$ be a self-adjoint operator with $\rho(B) = 1$ in Hilbert space $X$, not having $-1$ as its eigenvalue. Let be $Pf = 0$, where $P$ is an orthoprojection on the set of eigenvectors of operator $B$, corresponding to the eigenvalue $1$. Then residuals $x_n - Bx_n - f$ for successive approximations {\rm (\ref{102})}  at any starting condition $x_0 \in X$  converge to zero.}

We have to point out that condition $Pf = 0$ in this theorem is necessary, but in the general case not sufficient, for solving equation (\ref{101}). Consequently, residuals for successive approximations can converge to zero also in the case when the original equation has no solutions at all.

It follows from Theorem 1.2  that the convergence rate of residuals and corrections zero in this case is defined by the first residual properties $x_0 - Bx_0 - f$.

\vspace{0.3cm}

{\bf 1.3. Convergence оf errors, residuals and corrections in special subspaces.} As it is shown by simple examples and equations (\ref{106}), (\ref{108}) the convergence rate of successive approximations to the exact solution and that of residuals to zero considerably depends on initial approximation $x_0$ and  right-hand member $f$ of equation (\ref{101}). It is possible to estimate these convergence rates more exactly for functions $f$ from some (usually unclosed!) subspaces $\widetilde{X}$ of space $X$. Among such subspaces the simplest ones are the subspaces of <<sourcewise>> representable functions. These subspaces are defined with the help of some function $\theta(\lambda)$ singled out in the spectrum ${\rm sp}\, B$ оf operator $B$ as a set of elements $\theta(B)X$ of the type
 \begin{equation}\label{109}
 x = \theta(B)h \ \bigg(= \int\limits_{{\rm sp}\, B} \theta(\lambda) \, dE_\lambda h\bigg) \qquad (h \in X).
 \end{equation}
The set $\theta(B)X$ changes into a normed linear space, if the norm on its elements is defined as equation
 \begin{equation}\label{110}
 \|x\|_{\theta(B)X} = \inf\, \bigg\{\|h\|:\ h \in X, \ \theta(B)h = x\bigg\}.
 \end{equation}
It presents no difficulty to test that with this norm (but not with the initial one!) space $\theta(B)X$ is a Banach space.

Formula (\ref{106}) at $x_0 - x_* \in \theta(B)X$ is rewritten in the form
 \begin{equation}\label{111}
 \|x_n - x_*\|^2 =  \int\limits_{{\rm sp}\, B} |\lambda|^{2n} |\theta(\lambda)|^2 \, (dE_\lambda h,h).     \end{equation}
By virtue of the spectral theorem for self-adjoint operators, there follows from it inequality
 \begin{equation}\label{112}
 \|x_n - x_*\| \le \gamma_n \, \|x_0 - x_*\|_{\theta(B)X} \qquad (x_0 - x_* \in \theta(B)X),     \end{equation}
where
 \begin{equation}\label{113}
 \gamma_n = \max_{\lambda \in {\rm sp}\, B} \, |\lambda|^n |\theta(\lambda)|.
 \end{equation}

If $\gamma_n \to 0$ at $n \to \infty$, то (\ref{112}) gives the qualified estimate of the convergence rate of approximations (\ref{102}) for solving equation (\ref{101}) for all functions $x_0$ and $f$ at once, for which $x_0 - x_* \in \theta(B)X$. The latter condition is difficult to test, since $x_*$ is unknown. However, it is satisfied if $x_0 - Bx_0 - f \in \widetilde{\theta}(B)X$, where functions $\theta$ и $\widetilde{\theta}$ are connected by equality $\theta(\lambda) = (1 - \lambda)\widetilde{\theta}(\lambda)$. As a result, instead of  (\ref{112}) we have estimate
 \begin{equation}\label{114}
 \|x_n - x_*\| \le \widetilde{\gamma}_n \, \|x_0 - Bx_0 - f\|_{\widetilde{\theta}(B)X} \qquad (x_0 - Bx_0 - f \in \widetilde{\theta}(B)X),
 \end{equation}
where
 \begin{equation}\label{115}
 \widetilde{\gamma}_n = \max_{\lambda \in {\rm sp}\, B} \, |\lambda|^n |\widetilde{\theta}(\lambda)|.
 \end{equation}

Similarly, formula (\ref{108}) at $(I - B)x_0 - f \in \theta(B)X$  brings us to estimate
 \begin{equation}\label{116}
 \|x_n - Bx_n - f\| \le \gamma_n \, \|h\| \qquad (x_0 - Bx_0 - f = \theta(B)h, \ h \in X),
 \end{equation}
where sequence $(\gamma_n)$ is again defined by equality (\ref{113}).

The following holds

{\bf Theorem 1.3.} {\it Let $B$ be a self-adjoint operator with $\rho(B) = 1$ in Hilbert space $X$, not having $-1$ as its eigenvalue. If $\theta$  is a function with $\theta(\pm1) = 0$, determined on spectrum ${\rm sp}\, B$ then $\gamma_n \to 0$ and, сonsequently, at $x_0 - x_* \in \theta(B)X$ the convergence rate of approximations {\rm (\ref{102})} to consistent solution $x_*$ of equation {\rm (\ref{101})} is estimated by inequality {\rm (\ref{112})}. Further, if $\theta(\lambda) = (1 - \lambda)\widetilde{\theta}(\lambda)$ with $\widetilde{\theta}(\pm1) = 0$, то  $\widetilde{\gamma}_n \to 0$ and, hence, at $x_0 - Bx_0 - f \in X(\widetilde{\theta})$ the convergence rate of approximations {\rm (\ref{102})} to consistent solution $x_*$ of equation {\rm (\ref{101})} is defined by inequality} (\ref{114}).

{\bf Theorem 1.4.} {\it Let $B$ be a self-adjoint operator with $\rho(B) = 1$ in Hilbert space $X$, not having $-1$ as its eigenvalue. If $\theta$ is a function with $\theta(\pm1) = 0$, determined on spectrum ${\rm sp}\,  B$ then  $\gamma_n \to 0$ and, сonsequently, at $x_0 - Bx_0 - f \in \theta(B)X$ residual convergence rate for approximations {\rm (\ref{102})} to zero is estimated by inequality} (\ref{116}).

Both theorems follow from the following lemma.

{\bf Lemma 1.1.} {\it Let function $\vartheta(\lambda):\ [-1,1] \to {\Bbb R}$ satisfy the condition $\vartheta(\pm 1) = 0$. Then}
 \begin{equation*}
 \lim_{n \to \infty}\, \max_{-1 \le \lambda \le 1}\, |\lambda|^n |\vartheta(\lambda)| = 0.     \end{equation*}

$\square$ Let there be given $0 < \varepsilon < 1$. Then there exists such $\delta > 0$, that  at $1 - \delta < |\lambda| \le 1$ the inequality $|\vartheta(\lambda)| < \varepsilon$ is true. On the set $\{\lambda:\ |\lambda| \le 1 - \delta\}$ there holds inequality $|\lambda|^n |\vartheta(\lambda)| \le c(1 - \delta)^n$, where $c = \max\limits_{-1 \le \lambda \le 1}\, |\vartheta(\lambda)|$, and, hence, $|\lambda)|^n |\vartheta(\lambda)| < \varepsilon$ at $n > \D\frac{\ln (c^{-1}\varepsilon)}{\ln (1 - \delta)}$. But at $\lambda \in \{\lambda:\ 1 - \delta < |\lambda| \le 1\}$  inequality $|\lambda|^n |\vartheta(\lambda)| < \varepsilon$ also holds, and, consequently, this inequality holds at all $\lambda \in [-1,1]$. Since $\varepsilon$ is arbitrary, and $n$ does not depend on $\lambda$, then $|\lambda|^n |\vartheta(\lambda)| \to 0$ at $n \to \infty$ is uniform according to $\lambda \in [-1,1]$. $\blacksquare$

It should be noted that the conditions of theorems 1.3 и 1.4 соntain initial approximation $x_0$. If, as it is usually done, $x_0 = 0$, then the conditions of theorems 1.3 и 1.4 come down to assumptions concerning the solution itself $x_*$ or the given right-hand member $f$. The latter also holds when $x_0$ is taken nonzero, but <<good enough>> (in the examples, <<differentiable enough>>).

Finally, we consider that the assertions of theorems 1.3  и 1.4 essentially mean the convergence to zero according to the sequence norm of operators $B^n\theta(B)$ or the convergence to zero of the operator sequence $B^nT\theta(B)$, where $T$ is a quasi-inverse (possibly unlimited) operator for operator $(I - B)$ ($(I - B)T(I - B) = I - B$).

\vspace{0.3cm}

{\bf 1.4. Сonvergence in <<weakened>> norms.} In a number of jobs, while studying successive approximations, it is enough to determine their convergence in a weaker norm than the original norm of Hilbert space $X$.  Such norms can be exemplified by norm
 \begin{equation}\label{117}
 \|x\|_0 = \|Tx\|,
 \end{equation}
where $T$ is some noninvertible operator with $\ker T = 0$. Herewith, the simplest case is when operator $T$ is commutative with operator $B$ ($TB = BT$). Among such operators the simplest ones belong to the type
 \begin{equation}\label{118}
 T_\pi = \pi(B),
 \end{equation}
where $\pi(\lambda)$ is some bounded function, for which elements ${\rm Sp}\, B \cap \{\lambda:\ \pi(\lambda) = 0\}$ are not eigenvalues. In this case (\ref{117}) is the norm, because it follows from $Tx = 0$ that $x = 0$. The norms of such type are sometimes called weakened or relaxed generating ones. It should be noted that space $X$ with the norm (\ref{117}) is incomplete, if function $\pi^{-1}(\lambda)$ is unbounded on spectrum ${\rm Sp}\, B$.

We find it necessary to consider equality (\ref{105}):
 \begin{equation*}
 x_n - x_* = B^n(x_0 - x_*);
 \end{equation*}
resulting from (\ref{103}), (\ref{104}). Here $x_n$ are successive approximations,  $x_{n+1} = Bx_n + f$ с $x_0 \in X$ is the initial approximation to the solution of equation (\ref{101}), $x_*$ is the exact solution of equation (\ref{101}).

From this equation for the norm (\ref{117}) with $T$, defined by equation \ref{118} we have equation
 \begin{equation*}
 \|x_n - x_*\|_\pi  = \|\pi(B)B^n(x_0 - x_*)\|,
 \end{equation*}
and, further,
 \begin{equation*}
 \|x_n - x_*\|_0^2 = \int\limits_{{\rm Sp}\, B} |\pi(\lambda)|^2 |\lambda|^{2n} \, (dE_\lambda(x_0 - x_*),x_0 - x_*),
 \end{equation*}
from where,
 \begin{equation}\label{119}
 \|x_n - x_*\|_\pi \le \gamma_n \|x_0 - x_*\|,
 \end{equation}
where
 \begin{equation}\label{120}
 \gamma_n = \max_{\lambda \in {\rm Sp}\, B}\, |\pi(\lambda)| \, |\lambda|^n.
 \end{equation}

By using lemma 1.1 we come to the following assertion supplementing theorem 1.1.

{\bf Theorem 1.5.} {\it Let $B$ be a self-adjoint operator with $\rho(B) = 1$ in Hilbert space $X$, not having $-1$ as eigenvalue. Let $\pi(\pm 1) = 0$  and equation {\rm (\ref{101})} is solvable. Then, successive approximations {\rm (\ref{102})} at any initial condition $x_0 \in X$ converge in norm {\rm (\ref{117})} to solution $x_*$ of equation  {\rm (\ref{101})}, for which $Px_* = Px_0$, where $P$ is an orthoprojection on the set of eigenvectors of operator $B$, corresponding to eigenvalue $1$. Then, this convergence is uniform as regards $x_0 - x_* \in X$ on each bounded set.}

We underline that in the conditions of theorem 1.5, there is  no demand for the sourcewise representability of the exact solution or the right-hand member of equation (\ref{101}). We also note that in the conditions of theorem 1.5, the sequence of approximations (\ref{102}), in case of equation (\ref{101}) is not solved, can be fundamental in norm (\ref{117}). In other words, it can prove to be convergent in completion $X_\pi$ of space $X$ under the norm (\ref{117}), while this limit turns out to be the generalized solution of equation (\ref{101}).

Similarly to theorem 1.5, one proves the following theorem 1.6; at that, instead of the equation (\ref{105}) equation (\ref{108}) is used which also results from (\ref{103}), (\ref{104}):
 \begin{equation*}
 x_n - Bx_n - f = B^n(x_0 - Bx_0 - f);
 \end{equation*}
Here $x_n$ are successive approximations $x_{n+1} = Bx_n + f$ с $x_0 \in X$, $x_0$ is the initial approximation to the solution of equation (\ref{101}) (the solution itself may not exist at all).

{\bf Theorem 1.6.} {\it Let $B$ be a self-adjoint operator with $\rho(B) = 1$ in Hilbert space $X$, not having $-1$ as eigenvalue. Let $Pf = 0$, where $P$ is an orthoprojection on the set of eigenvectors of operator $B$, corresponding to the characteristic constant $1$. Then, residuals $x_n - Bx_n - f$ for successive approximations {\rm (\ref{102})} at any initial condition $x_0 \in X$  converge in the norm {\rm (\ref{117})} to zero. Consequently, this convergence is uniform in relation to $x_0 - Bx_0 - f \in X$ on each bounded set.}

\vspace{0.3cm}

{\bf 1.5. Сonvergence in errors of estimation.} Let now the conditions of theorem 1.1. be again satisfied for self-adjoint operator $B$. Let equation (\ref{101}) be solvable. In this case, the successive approximations (\ref{102}) сonverge to one of the solutions $x_*$ of equation (\ref{101}). Consider now, instead of exact successive approximations (\ref{102}) the approximations for the case when the right-hand member of the equation (\ref{101}) is set approximately, or when an error is made at every step of estimating these approximations. Both variants of such approximations are described well enough by equations
 \begin{equation}\label{121}
 \widetilde{x}_{n+1} = B\widetilde{x}_n + f_n \qquad (n = 0,1,2.\ldots)
 \end{equation}
assuming that $\|f_n - f\| \le \delta_n$ ($n = 0,1,2,\ldots$), where $(\delta_n)$ is some sequence of small positive numbers, bounded by the number $\delta$. From these equations and (\ref{102}), it directly follows
 \begin{equation*}
 \widetilde{x}_n = x_n + B^{n-1}(f_0 - f) + B^{n-2}(f_1 - f) \ldots + (f_{n-1} - f)     \end{equation*}
and, hence,
 \begin{equation*}
 \|\widetilde{x}_n - x_n\| \le \|B^{n-1}\| \, \|f_0 - f\| + \ldots + \|B\| \, \|f_{n-2} - f\| + \|f_{n-1} - f\| \le \delta_0 + \ldots + \delta_{n-2} + \delta_{n-1}.
 \end{equation*}
In this way,
 \begin{equation*}
 \|\widetilde{x}_n - x_*\| \le \|x_n - x_*\| + \|\widetilde{x}_n - x_n\|,
 \end{equation*}
and, consequently,
 \begin{equation}\label{122}
 \|\widetilde{x}_n - x_*\| \le \|x_n - x_*\| + (\delta_0 + \ldots + \delta_{n-2} + \delta_{n-1}),
 \end{equation}
where $x_*$ is the exact solution of equation (\ref{101}).

From inequalities (\ref{122}) the convergence $\widetilde{x}_n$ к $x_*$ does not follow, since the right-hand part in (\ref{122}) at $n \to \infty$  does not tend to zero (and, moreover, usually tends to infinity). However, in many cases it follows from these inequalities that, on the one hand, at quite big, but not too big, numbers $n$, the approximations (\ref{121}) come close enough to the exact solution $x_*$ of equation (\ref{101}). Moreover, these approximations for sequences $(\delta_n)$, sufficiently small in natural sense, <<fit>> the exact solution $x_*$ arbitrarily close!

In the conditions of theorem 1.1, at every initial approximation $x_0 \in X$ exact approximations $x_n$ сonverge to $x_*$, or, to put it differently, for some sequence of nonnegative numbers $\mu_n$ tending to zero, inequality
 \begin{equation*}
 \|x_n - x_*\| \le \mu_n.
 \end{equation*} holds.

We also recall that in the conditions of theorem 1.1, the sequence $(\mu_n)$ essentially depends on the initial condition $x_0 \in X$ and the right-hand member $f \in X$. However, theorem 1.3 also enables to describe some sets of initial conditions $x_0 \in X$ and right-hand members $f \in X$, for whose elements sequence $\mu_n$ can be chosen, independent of $x_0 \in X$ и $f \in X$.

Suppose
 \begin{equation*}
 \Delta_0 = 0, \qquad \Delta_n = \delta_0 + \ldots + \delta_{n-2} + \delta_{n-1} \quad (n = 1,2,\ldots).
 \end{equation*}
Then the inequality (\ref{122}) is transcribed as
 \begin{equation}\label{123}
 \|\widetilde{x}_n - x_*\| \le \mu_n + \Delta_n \qquad (n = 0,1,2,\ldots).
 \end{equation}

To estimate the <<smallness>> of sequence $(\delta_n)$, it is most convenient to suppose that the sequence $(\delta_n)$ belongs to some Banach space $L$ (with a monotonous, in the usual sense, norm) and to estimate this <<smallness>> by the norm $\|(\delta_n)\|_L$. It appears, the numbers $\delta_0 + \delta_1 + \ldots + \delta_{n-1}$ ($n = 1,2,\ldots$) can be regarded as meanings in sequence $(\delta_n)$ of linear functionals $\sigma_n$ ($n = 0,1,2,\ldots$), generated by sequence $(1,1,\ldots,1,0,\ldots)$, the first elements of which $n$ are equal to $1$, аnd the rest are equal to zero. According to the definition of the norms of functionals $\sigma_n$, the following inequalities
 \begin{equation}\label{124}
 \delta_0 + \delta_1 + \ldots + \delta_{n-1} \le \|\sigma_n\| \, \delta \qquad (\delta = \|(\delta_n)\|_L, \ n = 0,1,2,\ldots).
 \end{equation}
hold. From inequalities (\ref{123}) и (\ref{124}), there follow the estimates
 \begin{equation}\label{125}
 \|\widetilde{x}_n - x_*\| \le \mu_n + \|\sigma_n\| \, \delta \qquad (\delta = \|(\delta_n)\|, \ n = 0,1,2,\ldots)
 \end{equation}
The sequence $(\|\sigma_n)\|)$ is increasing; the examples below show that it can be both unbounded and bounded.

We can easily describe the behaviour peculiarities of sequence $(\mu_n + \|\sigma_n\| \delta)$ in the form of the following statement. Then it will be convenient for us to consider further on a more general sequence $(\mu_n + c \|\sigma_n\| \delta)$, where $c$  is some positive number.

{\bf Lemma 1.2.} {\it Let sequence $(\mu_n)$ tend to zero, while sequence $(\|\sigma_n\|)$ is nondecreasing. Then
 \begin{equation}\label{126}
 \lim_{n \to \infty, \, \|\sigma_n\| \delta \to 0}\ (\mu_n + c \|\sigma_n\| \delta) = 0.     \end{equation}
More exactly, let $\varepsilon > 0$ be given. Then there exists such $N(\varepsilon)$, that  at any $N_-$, $N_+$, for which $N(\varepsilon) \le N_- < N_+$ such $\delta(N_-,N_+)$ exists, that at $0 < \delta < \delta(N_-,N_+)$ the inequalities
 \begin{equation}\label{127}
 \mu_n + c \|\sigma_n\| \delta < \varepsilon \qquad n \in [N_-,N_+].
 \end{equation}
hold.}

In other words, at the given $\varepsilon > 0$ at sufficiently small $\delta > 0$  the inequality $\mu_n + \|\sigma_n\| \delta < \varepsilon$ is satisfied within arbitrarily distant and arbitrarily big change gaps $n$.

$\square$ The equality (\ref{126}) is obvious. Let now give $\varepsilon > 0$. To establish the inequality (\ref{127}), we first mention that at any $t$, $0 < t < \varepsilon$, at $n > N(t)$ the inequality $\mu_n < t$ is satisfied. Further, at the same $t$, we take arbitrary numbers $N_-$, $N_+$, for which $N(\varepsilon) \le N_- < N_+$ and then the number $\delta(N_-,N_+)$ so that at $n \in [N_-,N_+]$ the inequality
 \begin{equation*}
 \|\sigma_n\| < \frac{\varepsilon - t}{c\delta}.
 \end{equation*}
be satisfied. Then, at $\delta \le \delta(N_-,N_+)$ и $n \in [N_-,N_+]$
 \begin{equation*}
 \mu_n + \|\sigma_n\| \delta < t + \frac{\varepsilon - t}{c\delta} \cdot c\delta = \varepsilon. \qquad \qquad \blacksquare
 \end{equation*}

The relation  (\ref{126}) of lemma 1.2 is sometimes written in the form of
 \begin{equation}\label{128}
 \lim_{\delta \to 0}\, \min_{\nu \le n < \infty}\, \{\mu_n + c \|\sigma_n\| \delta\} = 0 \qquad (\nu \in {\Bbb N}).
 \end{equation}
However, without the additional assumption concerning the convergence of the sequence $(\mu_n)$ to zero, this relation is weaker than (\ref{126}).

We make another important remark. The inequalities (\ref{126}) turn out to be useful only in the cases, when at increasing $n$ the right-hand member $\mu_n + \|\sigma_n\| \delta$ decreases. The incident of decreasing the right-hand member in one step is equivalent to the inequality $\delta < \dfrac{\mu_n - \mu_{n+1}}{\|\sigma_{n+1}\| - \|\sigma_n\|}$. Тhus, the arguments considered show that the sequential computation of approximations (\ref{121}) prove to be  useful at $n \in [0,N]$ only if
 \begin{equation}\label{129}
 \delta < \frac{\mu_n - \mu_{n+1}}{c(\|\sigma_{n+1}\| - \|\sigma_n\|)} \qquad (n = 0,1,\ldots,N).     \end{equation}

In satisfying this correlation, one states that the correlating iteration method {\it quasi-converges}.

We once again point out that in case of quasi-convergence of iteration methods (\ref{121}) и (\ref{122}), one does not speak about the usual convergence of corresponding approximations to the exact solution. We can only assert that at sufficiently small $\delta$, these approximations happen to come close to the exact solution, and then, as a rule, move away from it; besides, the closer to the exact solution these approximations are, the less $\delta$ is. Moreover, if $\delta$ is not sufficiently small, then the use of approximations (\ref{121}) will turn out to be useless --- these approximations can move away from the exact solution.

It follows from the given considerations and lemma 1.2 that

{\bf Theorem 1.7.} {\it Let the conditions of the theorem {\rm 1.1} be satisfied and let the approximations {\rm (\ref{102})} be calculated with errors, not exceeding $\delta_n > 0$ at every step $n = 0,1,2,\ldots$, while $(\delta_n) \in L$, where $L$ is a Banach space of sequences with the monotonic norm. Then the approximations {\rm (\ref{121})} <<quasiсonverge>>, in the sense described above, to the corresponding solution $x_*$ of the equation {\rm (\ref{101})} (namely, the following relation holds}
 \begin{equation}\label{130}
 \lim_{n \to \infty, \, \|\sigma_n\| \delta \to 0}\ \|\widetilde{x}_n - x_*\| = 0.
 \end{equation}

We also observe that in the above-mentioned  <<paradoxical case>> $\mu_n = 0$,  it turns out that the initial approximation $x_0$ coincides with the solution $x_*$. It is in this case that the arguments about  the sequence $(\mu_n + n\delta)$, given above, degenerate, and the estimate (\ref{126}) becomes useless. However, it should be like that, if the initial approximation coincides with the exact solution $x_*$, it is useless to accurately define this approximation by any iteration procedures.

Now it remains to give the formulas for the norms $\|\sigma_n\|$ of the functionals $\sigma_n$ ($n = 0,1,\ldots$) for classical spaces mentioned above. For the spaces $\ell_p$ ($1 \le p \le \infty$) there happen to be equations
 \begin{equation*}
 \|\sigma_n\| = n^\frac1{p'} \qquad \bigg(\frac1p + \frac1{p'} = 1, \ n = 0,1,2,\ldots\bigg).
 \end{equation*}

One should observe two special cases in this equation when $p = \infty$ и $p = 1$. In the first one the condition $(\delta_n) \in \ell_\infty$ means that errors are made in calculations which do not exceed number $\delta = \|(\delta_n)\|_{\ell_\infty}$; in this case $\|\sigma_n\| = n$ for all $n = 0,1,2,\ldots$. We also observe that the assumption as to $\delta_n \to 0$ (or, otherwise, $(\delta_n) \in \ell_\infty^\circ = c_0 \subset \ell_\infty$) does not result in clarifying the behaviour of the norm sequence $(\sigma_n)$, both sequences of the norms for spaces $c_0$ и $\ell_\infty$ coincide. In the other case, when $(\delta_n) \in \ell_1$, the sequence of norms $(\|\sigma_n\|)$ turns out to be bounded!

For spaces ${\sf m}(\omega)$ of the sequences bounded by weight весом $\omega$ ($\omega = (\omega_0,\omega_1,\omega,\ldots)$, $\omega_k > 0$, $k = 0,1,2,\ldots$) the formulas
 \begin{equation*}
 \|\sigma_n\| = \sum_{k=0}^{n-1} \frac1{\omega_k} \qquad (n = 0,1,2,\ldots)
 \end{equation*}
hold. In the particular case, when $\omega = (1,2^\nu,\ldots,(k-1)^\nu,\ldots)$ the following equations
 \begin{equation*}
 \|\sigma_n\| = \sum_{k=0}^{n-1} \frac1{k^\nu} \qquad (n = 0,1,2,\ldots)
 \end{equation*}
hold. In this equation, one should also observe a special case when $\nu > 1$. In this case, the sequence of norms $(\|\sigma_n\|)$, as well as in the case of space $\ell_1$, also appears to be bounded: $\|\sigma_n\|_{{\sf m}(\omega)} \le \zeta(\nu)$ ($n = 0,1,2,\ldots$); here $\zeta(\cdot)$ is Riemann function.

\vspace{0.3cm}

{\bf 1.6. The main example.} We can take as an example in space $X = L_2(\Omega)$, where $\Omega$ is a closed set of the segment $[-1,1]$  c $1 \in \Omega$ (или $-1 \in \Omega$), the equation
 \begin{equation*}
 x(t) = tx(t) + f(t).
 \end{equation*}
This equation is solved in $X$, if and only if $(1 - t)^{-1}f(t) \in L_2(\Omega)$. The successive approximations (\ref{102}) in this case appear to be
 \begin{equation*}
 x_{n+1}(t) = tx_n(t) + f(t)
 \end{equation*}
or, which is just the same,
 \begin{equation*}
 x_n(t) = t^nx_0(t) +  (1 + t + t^2 + \ldots + t^{n-1})f(t)     .
 \end{equation*}
They converge in $X$ (at any $x_0(t) \in L_2(\Omega)$) to the function $(1 - t)^{-1}f(t)$, which is, under the assumption of the solvability the equation, belongs to $L_2(\Omega)$. The equation in this example is not correct. The similar situation takes place if $X = L_2(\Omega,\sigma)$, where $\sigma$ is  some measure on $\Omega$, while $\sigma(\{-1\}) = 0$.

The cited example has a sufficiently general character --- it is a known fact that every self-adjoint operator with a simple spectrum is similar to the operator of multiplying by an independent argument in the space $L_2(\Omega,\sigma)$ for the suitable choice of measure $\sigma$. For self-adjoint operators $B$ with the non-simple spectrum the similar assertion also holds, but one has to take here a topologically complex disjunctive union of segments $[-1,1]$ as $\Omega$.

\vspace{0.5cm}

\centerline{\bf \S\ 2. The equations of the first order}

\vspace{0.5cm}

{\bf 2.1.} {\bf The convergence principle.} Let $A$ be a self-adjoint operator in Hilbert space $X$. Let us consider the linear equation
 \begin{equation}\label{201}
 Ax = y,
 \end{equation}
where $y \in X$. We are interested in the case when $0$ is the point of spectrum ${\rm Sp}\, A$ of operator $A$.

Let $\phi(\lambda)$ be some real and analytical function on the spectrum of the operator $A$, which takes value $1$ at zero point; then
 \begin{equation*}
 \phi(\lambda) = 1 - \lambda \psi(\lambda),
 \end{equation*}
where $\psi(\lambda)$ is also a real and analytical function on ${\rm Sp}\, A $. Polynomials or rational functions can serve as the most obvious examples of such functions.

For each function $\phi(\lambda)$ of the type described above one can define the operator $\phi(A)$; it is also a self-adjoint one. Operator $\psi(A)$ is also defined. The equation
 \begin{equation*}
 x - \phi(A)x = \psi(A)Ax.
 \end{equation*}
is obvious. From this equation it follows that every solution $x$ of the equation (\ref{201}) is the solution to the equation
 \begin{equation}\label{202}
 x = \phi(A)x + \psi(A)y.
 \end{equation}
An inverse also holds, but under the additional assumption that $0$ is not the eigenvalue of the operator $\psi(A)$. Actually, (\ref{202}) can be rewritten as follows
 \begin{equation*}
 \psi(A)(Ax - y) = 0,
 \end{equation*}
from where it follows that $x$ is also the solution of the equation (\ref{201}). The assumption that $0$ is not the eigenvalue of operator $\psi(A)$ is equivalent to the one that $1$ is not the eigenvalue of operator $\phi(A)$. The latter, obviously, means the solution of equation (\ref{202}), if it exists, is unique. Thus, if the equation (\ref{201}) has the unique solution $x_*$, then it is the unique solution of equation (\ref{202}), and inversely, if the equation (\ref{202}) has the unique solution, then it will be the unique solution of the equation (\ref{201}). It should be observed that in the general case (without the assumption that $0$ is not the eigenvalue of operator $A$) in case of solvability of the equation (\ref{201}) the solution $x$ of the equation (\ref{202}) is not obligatorily the solution of the equation (\ref{201}), however, the solution of the equation (\ref{201}) in this case, is sure to be the element $x + \phi(A)(\xi - x)$, where $\xi$  is an arbitrary solution of the equation (\ref{201}).

Therefore, instead of analyzing the solvability properties of the equation (\ref{201}) one can consider the equation (\ref{202}). However, the latter equation has the form of $x = Bx + f$ with $B = \phi(A)$, $f = \psi(A)y$, and for its analysis one can naturally use the theorem of M.A. Krasnosel'skii mentioned above and all results from \S\ 1. The conditions of the latter will be satisfied if $\|\phi(A)\| = 1$ и $-1$ is not the eigenvalue of operator $\phi(A)$. Since, by virtue of Danford theorem, \cite{DanShv} ${\rm Sp}\, \phi(A) = \phi({\rm Sp}\, A)$, and the operator $\phi(A)$ is self-adjoint, then the equation $\|\phi(A)\| = 1$ is equivalent to the inequality
 \begin{equation}\label{203}
 |\phi(\lambda)| \le 1 \qquad (\lambda \in {\rm Sp}\, A)
 \end{equation}
(let us recall that $\phi(0) = 1$ and, therefore, (\ref{203}) means $\|\phi(A)\| = 1$). The second condition means that not a root of the equation $\phi(\lambda) + 1 = 0$ is the eigenvalue of operator $A$. Thus, the following holds

{\bf Theorem 2.1.} {\it Let  $A$ be a self-adjoint operator in Hilbert space $X$ and its range of values is not closed. Let $\phi(\lambda)$ be the analytical function of the environment ${\rm Sp}\, A$, for which

{\rm a)} $\phi(\lambda) = 1 - \lambda\psi(\lambda)$;

{\rm b)} $|\phi(\lambda)| \le 1$ \  ($\lambda \in {\rm Sp}\, A$);{\rm c)} zeroes of the function $\phi(\lambda) + 1$ are not the eigenvalues of the operator $A$.  Consequently, if the equation {\rm (\ref{201})} is solved, the successive approximations
 \begin{equation}\label{204}
 x_{n+1} = \phi(A)x_n + \psi(A)y \qquad (n = 0,1,2,\ldots)
 \end{equation}
converge to one of the solutions of equation} (\ref{201}).

There, naturally arises the question of the convergence rate of the approximations (\ref{204}). From theorem 1.1, it follows that in the general case this rate can be arbitrarily slow. To make the picture complete, we give here the calculations from \S\ 1, modified directly for the equation (\ref{202}). From (\ref{204}), it clearly follows
 \begin{equation}\label{205}
 x_n = \phi^n(A)x_0 + (E + \phi(A) + \phi^2(A) + \ldots + \phi^{n-1}(A))\psi(A)y \qquad (n = 0,1,2,\ldots),
 \end{equation}
and from (\ref{202})
 \begin{equation}\label{206}
 x_* = \phi^n(A)x_* + (E + \phi(A) + \phi^2(A) + \ldots + \phi^{n-1}(A))\psi(A)y \qquad (n = 0,1,2,\ldots),
 \end{equation}
Subtracting (\ref{206}) from (\ref{205}), we obtain
 \begin{equation}\label{207}
 x_n - x_* = \phi^n(A)(x_0 - x_*) \qquad (n = 0,1,2,\ldots)
 \end{equation}
and, further,
 \begin{equation}\label{208}
 \|x_n - x_*\|^2 = \int\limits_{{\rm Sp}\, A} |\phi(\lambda)|^{2n} \, (dE_\lambda(x_0 - x_*),x_0 - x_*).
 \end{equation}
From the formula (\ref{208}) the convergence of the approximations $x_n$ to $x_*$ follows, by virtue  of Lebesgue theorem of the limiting process under integral sign for the sequence almost always converging to zero. As it has been noted, it follows from this formula that this sequence can turn out to be arbitrarily slow and considerably depend on the the properties of <<smoothness>> of the initial error $x_0 - x_*$, while  the latter can depend on the properties of <<smoothness>> of the right-hand member $y$ and the <<incorrecftness>> properties of the operator $A$. However, it should also be noted that this convergence is the faster, the  <<less>> the function $\phi(\lambda)$ on the spectrum of the operator $A$ is.

\vspace{0.3cm}

{\bf 2.2.} {\bf The convergence of residuals and corrections.} We consider now the behaviour of residuals $Ax_n - y$ and corrections $x_{n+1} - x_n = \phi(A)x_n + \psi(A)y - x_n$ for approximations (\ref{204}).

From (\ref{205}) it follows
 \begin{equation*}
 Ax_n = \phi^n(A)Ax_0 + (E + \phi(A) + \phi^2(A) + \ldots + \phi^{n-1}(A))\psi(A)Ay =
 \end{equation*}
 \begin{equation*}
 = \phi^n(A)Ax_0 + (E + \phi(A) + \phi^2(A) + \ldots + \phi^{n-1}(A))(E - \phi(A))y =     \end{equation*}
 \begin{equation*}
 = \phi^n(A)Ax_0 + (E - \phi^n(A))y,
 \end{equation*}
and, hence,
 \begin{equation}\label{209}
 Ax_n - y = \phi^n(A)(Ax_0 - y).
 \end{equation}
It follows from this equation that
 \begin{equation}\label{210}
 \|Ax_n - y\|^2 = \int\limits_{{\rm Sp}\, A} |\phi(\lambda)|^{2n} \, (dE_\lambda(Ax_0 - y),Ax_0 - y). \end{equation}

Analogously, from (\ref{205}) for corrections $x_{n+1} - x_n$ we have
 \begin{equation*}
 x_{n+1} - x_n = \phi^n(A)(\phi(A) - E)x_0 + \phi^n(A)\psi(A)y = \phi^n(A)(\phi(A)x_0 + \psi(A)y - x_0)
 \end{equation*}
or
 \begin{equation}\label{211}
 x_{n+1} - x_n = \phi^n(A)(x_1 - x_0).
 \end{equation}
Consequently,
 \begin{equation}\label{212}
 \|x_{n+1} - x_n\|^2 = \int\limits_{{\rm Sp}\, A} |\phi(\lambda)|^{2n} \, (dE_\lambda(x_1 - x_0),x_1 - x_0).
 \end{equation}

Consequently, from (\ref{210}) and (\ref{212}) we obtain the following assertion:

{\bf Theorem 2.2.} {\it Let the conditions of theorem {\rm 2.1}be satisfied. Let $Py = 0$, where $P$ is an orthoprojection on the set of the eigenvectors of the operator $\phi(A)$, corresponding to the eigenvalue $1$. Then, the residuals $Ax_n - y$ and corrections $x_{n+1} - x_n$ for successive approximations {\rm (\ref{204})} at any intial condition $x_0 \in X$  сonverge to zero.}

Here one should also note that residuals and corrections converge to zero without the assumption of the solvability of the equation (\ref{201}).

\vspace{0.3cm}

{\bf 2.3.} {\bf Convergence on subspaces.} As it is shown by the equations (\ref{208}), (\ref{210}), (\ref{212}) the rate of convergence of the successive approximations (\ref{202}) to the exact solution of the equation (\ref{201}), сorrespondingly, the convergence rate of the residuals and corrections to zero greatly depends on the right-hand member $y$ of the equation (\ref{201}) and the initial condition $x_0$. However, the rates of these convergences can be specified, if the right-hand members $y$ of the equation and, consequently, the initial conditions $x_0$ are taken from some subspaces $\widetilde{X}$ of the space $X$. The simplest subspaces among this kind are the subspaces of the sourcewise representable functions mentioned above. Namely, we consider the cases when the right-hand members $y$ of the equation and, hence, the initial conditions $x_0$ lie in the subspaces $\theta(A)X$, which are defined by the operator $A$,  exactly in the same way as the spaces $\theta(B)X$ were defined with the help of some function $\theta(\lambda)$ determined on ${\rm Sp}\, A$ as a set of elements of the type
 \begin{equation*}
 x = \D\int\limits_{{\rm Sp}\, A} \theta(\lambda) \, dE_\lambda h \qquad (h \in X)
 \end{equation*} 
with the norm
 \begin{equation*}
 \|x\|_{\theta(A)X} = \inf\, \big\{\|h\|:\ h \in X, \ \theta(A)h = x\big\}.
 \end{equation*}
Similar to the spaces $\theta(B)X$, we assume for the spaces $\theta(A)X$ that the zeroes of the function $\theta(\lambda)$ are not eigenvectors of the operator $A$.

In the assumption that $y \in \theta(A)X$ exists $h \in X$, for which $x = \theta(A)h$. Then
 \begin{equation}\label{213}
 \|x_n - x_*\|^2 = \int\limits_{{\rm Sp}\, A} |\phi(\lambda)|^{2n} |\theta(\lambda)|^2 \,   (dE_\lambda h,h).
 \end{equation}
Hence,
 \begin{equation}\label{214}
 \|x_n - x_*\| \le \gamma_n \|x_0 - x_*\|_{\theta(A)X} \qquad (x_0 - x_* \in \theta(A)X),     \end{equation}
where $\gamma_n = \D\max\limits_{\lambda \in {\rm Sp}\, A}\, |\phi(\lambda)|^n |\theta(\lambda)|$.

If $\gamma_n \to 0$ at $n \to \infty$, then (\ref{214}) gives the qualified estimate of the convergence rate of approximations (\ref{204}) to the solution of the equation (\ref{201}) immediately for all functions $x_0$ и $y$, for which $x_0 - x_* \in \theta(A)X$. The condition $x_0 - x_* \in \theta(A)X$ is difficult to check on, since  $x_*$ is unknown. However, it is satisfied if $Ax_0 - y \in \widetilde{\theta}(A)$, where the functions $\theta$ and $\widetilde{\theta}$ are connected by the equation $\theta(\lambda) = \lambda \widetilde{\theta}(\lambda)$. In this case, consequently, instead of (\ref{214}) we have the estimate
 \begin{equation}\label{215}
 \|x_n - x_*\| \le \widetilde{\gamma}_n \|Ax_0 - y\|_{\widetilde{\theta}(A)X} \qquad (Ax_0 - y \in \widetilde{\theta}(A)X),
 \end{equation}
where $\widetilde{\gamma}_n = \D\max\limits_{\lambda \in {\rm Sp}\, A}\, |\phi(\lambda)|^n |\widetilde{\theta}(\lambda)|$.

It is natural that to prove $\gamma_n \to 0$ at $n \to \infty$ and $\widetilde{\gamma}_n \to 0$ at $n \to \infty$ we will need an analogue of lemma 1.1:

{\bf Lemma 2.1.} {\it Let the function $\phi(\lambda):\ {\rm Sp}\, A \to {\Bbb R}$ satisfy the conditions of the theorem {\rm 2.1} and the function $\theta(\lambda):\ {\rm Sp}\, A \to {\Bbb R}$ is such that from $|\phi(\lambda)| = 1$ it follows that $\theta(\lambda) = 0$. Then}
 \begin{equation*}
 \lim_{n \to \infty}\, \max_{\lambda \in  {\rm Sp}\, A}\, |\phi(\lambda)|^n |\theta(\lambda)| = 0.     \end{equation*}

The proof of this lemma is absolutely analogous to the proof of lemma 1.1.

By virtue of the things mentioned above and lemma 2.1, it follows

{\bf Theorem 2.3.} {\it Let the conditions of theorem {\rm 2.1} be satisfied. Then:

{\rm a)} if $\theta$ is the function defined on the spectrum ${\rm Sp}\, A$, for which from $|\phi(\lambda)| = 1$ there follows $\theta(\lambda) = 0$, то $\gamma_n \to 0$ and, consequently, at $x_0 - x_* \in \theta(A)X$ the convergence rate of the approximations {\rm (\ref{204})} to the corresponding solution $x_*$ of the equation {\rm (1)} is estimated by the inequality {\rm (\ref{214})};

{\rm b)} if $\theta$ is the function defined on the spectrum ${\rm Sp}\, A$ for which from $|\phi(\lambda)| = 1$ it follows $\widetilde{\theta}(\lambda) = 0$, где $\widetilde{\theta}(\lambda) = \lambda^{-1}\theta(\lambda)$, то $\widetilde{\gamma}_n \to 0$, hence, at $Ax_0 - y \in \widetilde{\theta}(A)X$ the convergence rate of the approximations {\rm (\ref{204})} to the corresponding solution $x_*$ of the equation {\rm (\ref{201})} is estimated by the inequality {\rm (\ref{215})}.}

The formulas (\ref{210}) and (\ref{212}), in their turn, result in the estimates
 \begin{equation}\label{216}
 \|Ax_n - y\| \le \gamma_n \|Ax_0 - y\|_{\theta(A)X} \qquad (Ax_0 - y \in \theta(A)X),     \end{equation}
 \begin{equation}\label{217}
 \|x_{n+1} - x_n\| \le \gamma_n \|x_1 - x_0\| \qquad (x_1 - x_0 \in \theta(A)X),
 \end{equation}
where the sequence $\gamma_n$ is again defined by the equation (\ref{214}). From these considerations and lemma 2.1 again it follows

{\bf Theorem 2.4.} {\it Let the conditions of the theorem {\rm 2.1} be satisfied and let $\theta$ be a function defined on the spectrum ${\rm Sp}\, A$, for which from $|\phi(\lambda)| = 1$ it follows $\theta(\lambda) = 0$. Then $\gamma_n \to 0$ and, consequently, at $Ax_0 - y \in \theta(A)X$ the convergence rate of the residuals for approximations {\rm (\ref{204})} to zero is estimated by the inequality {\rm (\ref{216})} and at $x_1 - x_0 \in \theta(A)X$ the convergence rate of the residuals for approximations {\rm (\ref{204})} to zero is estimated by the inequality} (\ref{217}).

\vspace{0.3cm}

{\bf 2.4.} {\bf Сonvergence in weakened norms.} We continue to study the behaviour of successive approximations $x_{n+1} = \phi(A)x_n + \psi(A)y$ for the linear operator equation $Ax = y$ with a self-adjoint operator $A$, active in Hilbert space $X$ in case when $0$ is a point of the spectrum of operator $A$. In a number of problems studying successive approximations it is sufficient  to determine their convergence in the norm which is weaker than the initial norm of Hilbert space $X$.  Similarly to what has been done, we shall consider the norms
 \begin{equation}\label{218}
 \|x\|_0 = \|Tx\|
 \end{equation}
where $T$ is an  operator with ${\rm Ker}\, T = 0$ and such that $TA = AT$. In repeating considerations п. 1.4, we will restrict ourselves by operators of the type 
 \begin{equation}\label{219}
 T = \pi(A)
 \end{equation}
where $\pi$ is a function, positive on ${\rm Sp}\, A$, the zeroes of which are not the eigenvalues of operator $A$. In this case (\ref{218}) is the norm, as from $Tx = 0$ it obviously follows that $x = 0$.

Let us refer to [1], that the equations
 \begin{equation*}
 x_n = \phi^n(A)x_0 + (E + \phi(A) + \phi^2(A) + \phi^{n-1}(A))\psi(A)y,
 \end{equation*}
 \begin{equation*}
 x_* = \phi^n(A)x_* + (E + \phi(A) + \phi^2(A) + \phi^{n-1}(A))\psi(A)y.
 \end{equation*}
hold. Thus
 \begin{equation}\label{220}
 x_n - x_* = \phi^n(A)(x_0 - x_*) \qquad (n = 0,1,2,\ldots).
 \end{equation}
hold. Here $x_n$ are the successive approximations $x_{n+1} = \phi(A)x_n + \psi(A)y$ с $x_0 \in X$, $x_*$ is the exact solution of the equation $Ax = y$.

From the equation (\ref{220}) for norm (\ref{218}) (with $T$, defined by the equation (\ref{220})) there follows  the equation
 \begin{equation*}
 \|x_n - x_*\|_{\pi(A)X} = \|\pi(A)\phi^n(A)(x_0 - x_*)\|,
 \end{equation*}
and, further,
 \begin{equation*}
 \|x_n - x_*\|_{\pi(A)X}^2 = \int_{{\rm Sp}\, A} |\pi(\lambda)|^2 |\phi(\lambda)|^{2n} \, (dE_\lambda(x_0 - x_*),x_0 - x_*),
 \end{equation*}
hence,
 \begin{equation*}
 \|x_n - x_*\|_{\pi(A)X} \le \gamma_n \|x_0 - x_*\|,
 \end{equation*}
where $\gamma_n = \max_{\lambda \in {\rm Sp}\, A} \, |\pi(\lambda)| |\phi(\lambda)|^n$.

By repeating the statements from п. 1.4, we arrive at the following assertion supplementing the theorem of M.A. Krasnosel'skii.

{\bf Theorem 2.5.} {\it Let the conditions of theorem {\rm 2.1} be satisfied. Let $\pi(\pm 1) = 0$ and the equation $Ax = y$ be solved. Then the successive approximations $x_{n+1} = \phi(A)x_n + \psi(A)y$ at any initial condition $x_0 \in X$ сonverge in norm {\rm (\ref{218})} to the solution $x_*$ to the equation $Ax = y$, for which $Px_* = Px_0$, where $P$ is the orthoprojection on the set of eigenvectors of the operator $A$, corresponding to the eigenvalue $0$. Then this convergence is uniform with respect to $x_0 - x_* \in X$ from every bounded set.}

It is sufficient to show that $\gamma_n \to 0$ at $n \to \infty$. But this fact follows immediately from lemma 2.1, in which the function  $\theta$ should be changed for the function $\pi$.

We underline that under the conditions of theorem 2.5 the requirement for the sourcewise representability of the exact solution or the right-hand member of the equation (\ref{201}) is missing.

Theorem 2.5 is the analogue of theorem 1.5. The constructions described above allow to formulate the analogues of theorem 1.6 оn the convergence to zero of the residuals and corrections in norms (\ref{218}) at the corresponding choice of functions $\pi$ for the equations of the first order (\ref{201}).  We confine ourselves here only to the corresponding  definition.

{\bf Theorem 2.6.} {\it Let the conditions of theorem {\rm 2.1} be satisfied. Let  $\pi(\pm 1) = 0$ and  $Py = 0$, where $P$ is an orthoprojection on the subspace of the eigenvectors of the operator  $A$, corresponding to the eigenvalue $0$. Then the residuals $Ax_n - y$ and the corrections $x_{n+1} - x_n$ for successive approximations {\rm (\ref{204})} at any initial condition $x_0 \in X$  converge in norm {\rm (\ref{219})} to zero. Also this convergence is uniform with respect to $Ax_0 - y \in X$ and, consequently, $x_1 - x_0$ on every bounded set.}

\vspace{0.3cm}

{\bf 2.5.} {\bf The convergence of approximations at imperfect data and and in the occurrence of errors.} Let the conditions of theorem 2.1 be again satisfied for the self-adjoint operator $A$, while $\|\phi(A)\| = 1$ and, consequently, $\rho(\phi(A)) = 1$. Let the equation (\ref{201}) be solved. In this case the successive approximations (\ref{204}) converge to one of the solutions $x_*$ to the equation (\ref{201}). Instead of exact approximations (\ref{204}), we consider now the approximations for the case when the right-hand member of the equation (\ref{201}) is calculated at every step $n$ with an error not exceeding $\delta_n$. These new approximations $\widetilde{x}_n$ are written in the form
 \begin{equation}\label{221}
 \widetilde{x}_{n+1} = \phi(A)\widetilde{x}_n + \psi(A)y_n \qquad (n = 0,1,2,,\ldots)
 \end{equation}
with  the approximate right-hand member $y_n$, $\|y_n - y\| \le \delta_n$. As it is easily seen, from the equalities (\ref{204}), there follow the equalities
 \begin{equation*}
 \widetilde{x}_n  = \phi^n(A)x_0 + (\psi(A)y_{n-1} + \phi(A)\psi(A)y_{n-2} + \ldots + \phi^{n-1}\psi(A)y_0)
 \end{equation*}
valid at all $n = 0,1,2,\ldots)$.

Consequently from (\ref{204})
 \begin{equation*}
 \widetilde{x}_n - x_n = \psi(A)(y_{n-1} - y) + \phi(A)\psi(A)(y_{n-2} - y) + \ldots + \phi^{n-1}\psi(A)(y_0 - y),
 \end{equation*}
and, by virtue of $\|\phi(A)\| = 1$,
 \begin{equation*}
 \|\widetilde{x}_n - x_n\| \le \|\psi(A)\| \, \|y_{n-1} - y\| + \|\psi(A)\| \, \|y_{n-2} - y\| + \ldots + \|\psi(A)\| \, \|y_0 - y\|,
 \end{equation*}
and, finally,
 \begin{equation*}
 \|\widetilde{x}_n - x_n\| = \|\psi(A)\| \, (\delta_0 + \delta_1 + \ldots + \delta_{n-1}) \qquad (n = 0,1,2,\ldots).
 \end{equation*}

Since $\|\widetilde{x}_n - x_*\| \le \|x_n - x_*\| + \|\widetilde{x}_n - x_n\|$, then it follows from the last inequality
 \begin{equation}\label{222}
 \|\widetilde{x}_n - x_*\| = \|x_n - x_*\| + \|\psi(A)\| \, (\delta_0 + \delta_1 + \ldots + \delta_{n-1}) \qquad (n = 0,1,2,\ldots).
 \end{equation}

Let
 \begin{equation}\label{223}
 c = \max_{\lambda \in {\rm Sp}\, A}\, |\psi(\lambda)|.
 \end{equation}
 It follows from the spectral theorem for self-adjoint operators that this number coincides with $\|\psi(A)\|$. Therefore, from (\ref{222}) и (\ref{223}) the analogous (\ref{125}) estimate     \begin{equation}\label{224}
 \|\widetilde{x}_n - x_*\| = \|x_n - x_*\| + c(\delta_0 + \delta_1 + \ldots + \delta_{n-1}) \qquad (n = 0,1,2,\ldots)
 \end{equation}
follows. One can apply lemma 1.2. to the right-hand member of this inequality. From it correlations (\ref{126}) and (\ref{128}) follow and, further, the analogue of theorem 1.7. In other words, the following holds

{\bf Theorem 2.7.} {\it  Let the conditions of theorem {\rm 2.1} be satisfied, and let the approximations {\rm (\ref{204})} at every step $n = 0,1,2,\ldots$ are calculated with errors not exceeding $\delta_n > 0$, while $(\delta_n) \in L$, where $L$ is a Banach space of sequences with monotonic norm.  Then, the approximations {\rm (\ref{221})} <<quasiconverge>> in the above sense to the corresponding solution $x_*$ to the equation {\rm (\ref{201})}, that is, the relation}
 \begin{equation}\label{224}
 \lim_{n \to \infty, \, \|\sigma_n\| \delta \to 0}\ \|\widetilde{x}_n - x_*\| = 0 holds.
 \end{equation}

\vspace{0.3cm}

{\bf 2.6.} {\bf The main example.} We can also consider here the equation
\begin{equation*}
 tx(t) = y(t).
 \end{equation*}
as an example in the space $X = L_2(\Omega)$, where $\Omega$ is some bounded closed set the straight line ${\Bbb R}$  c $0 \in \Omega$. This equation is solvable in $X$, if and only if $t^{-1}y(t) \in L_2(\Omega)$. The successive approximations (\ref{221}) in this case are such that
 \begin{equation*}
 x_{n+1}(t) = \phi(t)x_n(t) + \psi(t)y(t)
 \end{equation*}
or, which is just the same,
 \begin{equation*}
 x_n(t) = \phi(t)^nx_0(t) +  (1 + \phi(t) + \phi^2(t) + \ldots + \phi^{n-1}(t))\psi(t)y(t).
 \end{equation*}
On satisfying the conditions of the corresponding theorem of this paragraph, these successive approximations converge in $X$ (at any $x_0(t) \in L_2(\Omega)$) to the function $t^{-1}y(t)$, which, under the assumption of the solvability of the equation, belongs to $L_2(\Omega)$. The equation in this example is not correct. Similar to the equations of the second order, the analogous situation also takes place if $X = L_2(\Omega,\sigma)$, where $\sigma$ is some measure on $\Omega$, when $\sigma(\{-1\}) = 0$.

As it is mentioned in п. 1.6, the given example is of a sufficiently general character.

\vspace{0.5cm}\centerline{\bf \S\ 3. Partial iteration methods for the equations of the first order}

\vspace{0.5cm}

{\bf 3.1.} {\bf Implicit iterative schemes.} By choosing various functions $\phi(\lambda)$ и $\psi(\lambda)$, which satisfy conditions a), b), c) of the theorem 2.1, we obtain various iterative schemes of  approximate constructing the solutions of the equation (\ref{201}). We will confine ourselves here to several examples (cf. \cite{VaiVer,BakGon,SavMat,Mat}, where they are studied from the another aspect).

First of all, we consider (see [14]) the iteration method (\ref{204}), соrresponding to the polynomial
 \begin{equation}\label{301}
 \phi(\lambda) = (1 - \alpha\lambda)^k
 \end{equation}
($k$ is a natural number, $\alpha > 0$). For it
 \begin{equation}\label{302}
 \psi(\lambda) = \D\frac{1 - (1 - \alpha\lambda)^k}{\lambda}
 \end{equation}
and, further, condition b) of the theorem 2.1 is satisfied if ${\rm Sp}\, A \subseteq \bigg[0,\D\frac2\alpha\bigg]$, while condition c), if $\lambda = \D\frac2\alpha$ is not the eigenvalue of the operator $A$. Correspondingly, the iterations (\ref{204}) take the form
 \begin{equation}\label{303}
 x_{n+1} = (E - \alpha A)^kx_n + A^{-1}[E - (E - \alpha A)^k]y \qquad (n = 0,1,2,\ldots).
 \end{equation}

For this method it is most convenient to take function $\theta(\lambda) = \lambda^s$  as the function $\theta(\lambda)$ ($s$ is some positive number). This function satisfies the conditions of the theorem 2.3, if ${\rm Sp}\, A \subseteq [0,M]$, where $M < \D\frac2\alpha$. The sequence $(\gamma_n)$, defined by the equation $\gamma_n = \D\max\limits_{\lambda \in {\rm Sp}\, A}\, |\phi(\lambda)|^n |\theta(\lambda)|$, will be estimated at $M \le \D\frac1\alpha$ by the equality
 \begin{equation*}
 \gamma_n = \bigg(\frac{s}{\alpha(s + kn)}\bigg)^s \bigg(\frac{kn}{s + kn}\bigg)^{kn} \qquad (n = 0,1,2,\ldots)
 \end{equation*}
and at $\D\frac1\alpha < M < \D\frac2\alpha$ by the equation
 \begin{equation*}
 \gamma_n = \max\, \bigg\{\bigg(\frac{s}{\alpha(s + kn)}\bigg)^s \bigg(\frac{kn}{s + kn}\bigg)^{kn}, \, M^s(1 - \alpha M)^{kn}\bigg\} \qquad (n = 0,1,2,\ldots).
 \end{equation*}
As it is clearly seen, in both cases,the correlation
 \begin{equation}\label{304}
 \gamma_n \sim \bigg(\frac{s}{e\alpha k}\bigg)^s \, \frac1{n^s}.
 \end{equation} holds

We consider now (see [14]) the iteration method (\ref{204}), соrrelating tо the polynomial
 \begin{equation}\label{305}
 \phi(\lambda) = (1 - \alpha\lambda^k)
 \end{equation}
($k$ is a natural number, $\alpha > 0$). For it
 \begin{equation}\label{306}
 \psi(\lambda) =  \alpha\lambda^{k-1}.
 \end{equation}
Condition b) of the theorem 2.1 at even $k$  is satisfied if ${\rm Sp}\, A \subseteq \bigg[-\bigg(\D\frac2\alpha\bigg)^\frac1k,\bigg(\D\frac2\alpha\bigg)^\frac1k\bigg]$ and at odd $k$, if ${\rm Sp}\, A \subseteq \bigg[0,\bigg(\D\frac2\alpha\bigg)^\frac1k\bigg]$.  Further, condition c) is satisfied if $\lambda = \pm \bigg(\D\frac2\alpha\bigg)^\frac1k$ in the first case and $\lambda = \bigg(\D\frac2\alpha\bigg)^\frac1k$ in the second case are not the eigenvalue of operator $A$. Соrrespondingly, the iterations (\ref{204}) are as follows
 \begin{equation}\label{307}
 x_{n+1} = (E - \alpha A^k)x_n + \alpha A^{k-1}y \qquad (n = 0,1,2,\ldots).
 \end{equation}

For this method it is also most convenient to take the function $\theta(\lambda) = \lambda^s$ as the function $\theta(\lambda)$ ($s$ is some positive number). The conditions of the theorem 2.3 are satisfied if ${\rm Sp}\, A \subseteq [-M,M]$ at even $k$ and ${\rm Sp}\, A \subseteq [0,M]$ at odd $k$, where $M < \bigg(\D\frac2\alpha\bigg)^\frac1k$. The sequence $(\gamma_n)$, defined by the equation $\gamma_n = \D\max\limits_{\lambda \in {\rm Sp}\, A}\, |\phi(\lambda)|^n |\theta(\lambda)|$, will be estimated at $ M \le \bigg(\D\frac2\alpha\bigg)^\frac1k$ by the equation
 \begin{equation*}
 \gamma_n = \bigg(\frac{s}{\alpha(s + kn)}\bigg)^\frac{s}k \bigg(\frac{kn}{s + kn}\bigg)^n \qquad (n = 0,1,2,\ldots)
 \end{equation*}
and at $\bigg(\D\frac1\alpha\bigg)^\frac1k < M < \bigg(\D\frac2\alpha\bigg)^\frac1k$ by the equation
 \begin{equation*}
 \gamma_n = \max\, \bigg\{\bigg(\frac{s}{\alpha(s + kn)}\bigg)^\frac{s}k \bigg(\frac{kn}{s + kn}\bigg)^n, \, M^s(1 - \alpha M^k)^n\bigg\} \qquad (n = 0,1,2,\ldots).
 \end{equation*}
It is clearly seen that in both cases the correlation
 \begin{equation}\label{308}
 \gamma_n \sim \bigg(\frac{s}{e\alpha k}\bigg)^\frac{s}k \, \frac1{n^\frac{s}k}
 \end{equation}
holds.

The comparison of equations (\ref{304}) и (\ref{308}) shows  that, что at $\theta(\lambda) = \lambda^s$ method (\ref{303}) converges asymptotically better than method (\ref{307}).

{\bf 3.2.} {\bf The implicit iterative schemes.} Let us first consider the case $\phi(\lambda) = \D\frac1{1 + \alpha\lambda^k}$ ($\alpha > 0$) and, соrrespondingly, $\psi(\lambda) = \D\frac{\alpha\lambda^{k-1}}{1 + \alpha\lambda^k}$. In this case we deal with the implicit method of iterations defined by the formulas
 \begin{equation}\label{309}
 (E + \alpha A^k)x_{n+1} = x_n + \alpha A^{k-1}y \qquad (n = 0,1,2,\ldots).
 \end{equation}
Condition b) of the theorem 2.1, at even $k$, is satisfied at ${\rm Sp}\, A \subseteq (-\infty,\infty)$ (т.е., всегда) and at odd $k$, if ${\rm Sp}\, A \subseteq [0,\infty)$.  Further, condition c) is always satisfied.

In order to use the theorem 2.3, we again consider the case $\theta(\lambda) = \lambda^s$; here $s$ is any positive number if ${\rm Sp}\, A \subseteq [0,\infty)$ and  a rational positive number with an even denominator, if ${\rm Sp}\, A \cap (-\infty,0)$. Simple calculations show that $\gamma_n = 1$ for $n \le \D\frac{s}k$ and for $n > \D\frac{s}k$
 \begin{equation}\label{310}
 \gamma_n = \bigg(\frac{s}{nk - s}\bigg)^s\bigg(\frac{nk - s}{nk}\bigg)^n \sim \bigg(\frac{s}{e\alpha k}\bigg)^\frac{s}k \, \frac1{n^\frac{s}k}.
 \end{equation}

Similarly, we examine the case $\phi(\lambda) = \D\frac{1 - \alpha\lambda^k}{1 + \alpha\lambda^k}$ ($\alpha$ is a positive number) and, correspondingly, $\psi(\lambda) = \D\frac{2\alpha\lambda^{k-1}}{1 + \alpha\lambda^k}$. The iteration method (\ref{221}) in this case coincides with the implicit method of iterations defined by the equations
 \begin{equation}\label{311}
 (E + \alpha A^k)x_{n+1} = (E - \alpha A^k)x_n + 2\alpha A^{k-1}y \qquad (n = 0,1,2,\ldots).
 \end{equation}
Condition b) of the theorem 2.1 at even $k$ is satisfied at ${\rm Sp}\, A \subseteq (-\infty,\infty)$ (i.e., always) and at odd $k$, if ${\rm Sp}\, A \subseteq [0,\infty)$.  Further, condition c) is always satisfied.

For applying the theorem 2.3 we again consider the case $\theta(\lambda) = \lambda^s$; here $s$ is any positive number, if ${\rm Sp}\, A \subseteq [0,\infty)$ и and a rational positive number with  an even denominator, if ${\rm Sp}\, A \cap (-\infty,0) \ne \emptyset$. The calculations of constant $\gamma_n$ result in rather cumbersome formulas. Therefore, we confine ourselves to clarifying their asymptotic behaviour at $n \to \infty$. Actually, $\gamma_n = \max\limits_{\lambda \in {\rm Sp}\, A} \ \lambda^s(1 - \alpha\lambda^k)^n(1 + \alpha\lambda^k)^{-n}$. The derivative of the function $\xi(\lambda) = \lambda^s(1 - \alpha\lambda^k)^n(1 + \alpha\lambda^k)^{-n}$ is defined by the equation
 \begin{equation*}
 \xi'(\lambda) = \lambda^{s-1}(1 - \alpha\lambda^k)^{n-1}(1 + \alpha\lambda^k)^{-n-1} \, ((s - s\alpha^2\lambda^{2k}) - 2nk\alpha\lambda^k).
 \end{equation*}
At big $n$ this derivative turns into zero in the point $\lambda = \lambda_n$, для которой
 \begin{equation*}
 s(\alpha\lambda^k)^2 + 2nk(\alpha\lambda^k) - s = 0,
 \end{equation*}
from where 
 \begin{equation}\label{301}
 \alpha\lambda^k = \sqrt{\bigg(\frac{nk}{s}\bigg)^2 + 1} - \frac{nk}{s} = \frac1{\sqrt{\bigg(\dfrac{nk}{s}\bigg)^2 + 1} + \dfrac{nk}{s}} \sim \frac{s}{2k} \cdot \frac1n.
 \end{equation}
Hence, at $\lambda = \lambda_n$,
 \begin{equation}\label{312}
 \gamma_n = \lambda^s \bigg(\frac{1 - \alpha\lambda^k}{1 + \alpha\lambda^k}\bigg)^n \sim \bigg(\frac{s}{2e\alpha k}\bigg)^\frac{s}k \cdot \frac1{n^\frac{s}k}.
 \end{equation}

Finally, let us consider another case when $\phi(\lambda) = \D\frac{(1 - \alpha\lambda^k)^2}{1 + \alpha^2\lambda^{2k}}$, $\alpha > 0$, and, consequently, $\psi(\lambda) = \D\frac{2\alpha\lambda^{k-1}}{1 + \alpha^2\lambda^{2k}}$. In this case we obtain the following iteration method
 \begin{equation}\label{313}
 (E + \alpha^2 A^{2k})x_{n+1} = (E - \alpha A^k)^2x_n + 2\alpha A^{k-1}y \qquad (n = 0,1,2,\ldots).
 \end{equation}
At odd $k$, we obtain that the conditions of the theorem 2.1 are satisfied if ${\rm Sp}\, A \subseteq [0,\infty)$, while at even $k$ we get that the conditions of the theorem 2.1 are satisfied at all times.

For applying the theorem 2.3 we again consider the case $\theta(\lambda) = \lambda^s$; here $s$ is any positive number. The calculations of constant $\gamma_n$ сomes down here to the analysis of the roots of some cubic equation. However, the asymptotic behaviour of these roots is defined rather simply; it turns out that for $\lambda = \lambda_n$ the correlation
 \begin{equation*}
 \lambda \sim \bigg(\frac{s}{2k\alpha}\bigg)^\frac1k \cdot \frac1{n^\frac1k}
 \end{equation*}
holds, and therefore
 \begin{equation}\label{314}
 \gamma_n \sim \bigg(\frac{s}{2ke\alpha}\bigg)^\frac{s}k \cdot \frac1{n^\frac{s}k}.
 \end{equation}

The comparison of correlations (\ref{310}), (\ref{312}), (\ref{314}) shows that the convergence rate of all the three methods considered in this paragraph is asymptotically equal.

\vspace{0.5cm}\def\refname{\centerline{\normalsize \bf RERERENCES}}

\end{document}